\documentclass[a4paper,12pt]{amsart}
\usepackage{amssymb}
\usepackage{amsmath}
\usepackage[pdftex]{hyperref}

\newcommand{\hide}[1]{}
\newcommand{\ignore}[1]{}

\newcommand{\C}{\mathbb{C}}

\newcommand{\Q}{\mathbb{Q}}

\newcommand{\Legendre}[2]{\genfrac{(}{)}{}{}{#1}{#2}}

\DeclareMathOperator{\arcsinh}{arcsinh}

\DeclareMathOperator{\ord}{ord}

\newtheorem{dummy}{Dummy}

\newtheorem{lemma}[dummy]{Lemma}
\newtheorem{theorem}[dummy]{Theorem}
\newtheorem{prop}[dummy]{Proposition}
\newtheorem{cor}[dummy]{Corollary}

\theoremstyle{definition}

\theoremstyle{remark}
\newtheorem{rem}[dummy]{Remark}

\begin{document}
\bibliographystyle{amsalpha}

\author{S.~Mattarei}
\email{mattarei@science.unitn.it}
\urladdr{http://www-math.science.unitn.it/\~{ }mattarei/}
\address{Dipartimento di Matematica\\
  Universit\`a degli Studi di Trento\\
  via Sommarive 14\\
  I-38050 Povo (Trento)\\
  Italy}

\title{The sum $\sum_{k=0}^{q-1}\binom{2k}{k}$ for $q$ a power of $3$}

\begin{abstract}
We prove that
$\sum_{k=0}^{q-1}\binom{2k}{k}
\equiv
q^2\pmod{3q^2}$
if $q>1$ is a power of $3$,
as recently conjectured by Z.W.~Sun and R.~Tauraso.
Our more precise result actually implies that the value of
$(1/q^2)\sum_{k=0}^{q-1}\binom{2k}{k}$
modulo a fixed arbitrary power of $3$
is independent of $q$, for $q$ a power of $3$ large enough, and shows how such value can be efficiently computed.
\end{abstract}



\maketitle

\thispagestyle{empty}

\section{Introduction}

The sum of central binomial coefficients $\sum_{k=0}^{q-1}\binom{2k}{k}$, where $q$ is a power of a prime $p$,
has recently received some attention.
Because of the connection with the combinatorially more interesting Catalan numbers
$C_n=\binom{2n}{n}-\binom{2n}{n+1}$,
its variation $\sum_{k=0}^{q-1}\binom{2k}{k+d}$ has also been considered, where $d$ is a fixed integer.
In particular, the congruence
\begin{equation}\label{eq:PanSun}
\sum_{k=0}^{q-1}\binom{2k}{k+d}\equiv\Legendre{q-d}{3}\pmod{p}
\end{equation}
was first noted and proved
by H.~Pan and Z.W.~Sun in~\cite{PanSun},
in the special case where $q=p$.
Here $\Legendre{a}{3}$ is a Legendre symbol, hence uniquely determined by
$\Legendre{a}{3}\in\{0,\pm 1\}$ and
$\Legendre{a}{3}\equiv a\pmod{p}$.

Congruence~\eqref{eq:PanSun} has been extended in~\cite{SunTau:Catalan} to a congruence modulo $p^2$,
which generally involves additional terms.
However, in the special case $d=0$ (as well as the case $d=1$ if $p\neq 3$), congruence~\eqref{eq:PanSun} remains true modulo $p^2$ in the form stated above.
In particular, when $p=3$ the sum
$\sum_{k=0}^{q-1}\binom{2k}{k}$
is a multiple of $9$.
In fact, Z.W.~Sun and R.~Tauraso have noted that it is actually a multiple of $q^2$:
they conjecture in~\cite{SunTau:Catalan} that
\begin{equation}\label{eq:conj}
\frac{1}{q^2}
\sum_{k=0}^{q-1}\binom{2k}{k}
\equiv
1\pmod{3}
\end{equation}
if $q>1$ is a power of $3$.
In this note we confirm this conjecture by proving the following stronger result.
We write $H_{k-1}=\sum_{0<i<k}1/i$ for the harmonic numbers.

\begin{theorem}\label{thm:main}
Let $q=3^f$, and let $e<f$ be an integer such that $4\cdot 3^{e-1}+e\ge 2f+2$.
Then
\[
\frac{1}{q^2}
\sum_{k=0}^{q-1}\binom{2k}{k}
\equiv
-\beta\binom{2q}{q}
\pmod{3^{f-2e+2}},
\]
where $\beta$ is the $3$-adic integer
\[
\beta=\sum_{k=2}^{\infty}\frac{3^{k-1}}{k}
\binom{2k}{k}^{-1}
H_{k-1}
=
1+3+2\cdot 3^3+2\cdot 3^4+2\cdot 3^5+3^7+2\cdot 3^{11}+
\cdots.
\]
\end{theorem}

The condition $4\cdot 3^{e-1}+e\ge 2f+2$ is quite close to optimal for the validity of Theorem~\ref{thm:main}.
It can be slightly weakened at the expense of a more tedious argument, as we explain in Remark~\ref{rem:p-adic}.

The following simple proposition gives an upper bound of the number of terms required for the $3$-adic evaluation of $\beta$
with the precision needed in Theorem~\ref{thm:main},
depending on $f$ and $e$.
In particular, it shows that the number of necessary terms is $O(\log q)$.
However, a number of terms independent of $q$ will suffice if we choose $e$ at a fixed distance from $f/2$, as we will illustrate below.

\begin{prop}\label{prop:main}
With notation as in Theorem~\ref{thm:main} we have
\[
\beta\equiv
\sum_{k=2}^{s}\frac{3^{k-1}}{k}
\binom{2k}{k}^{-1}
H_{k-1}
\pmod{3^{f-2e+2}},
\]
where $s\ge 2$ is any integer such that
$3^s/(2s^2)\ge 3^{f-2e+1}$.
\end{prop}

Theorem~\ref{thm:main} provides the most precise $3$-adic evaluation of our binomial sum when $e$ is chosen as small as possible,
which means $e=\bigl\lfloor\log\bigl((f+1)/2\bigr)/\log 3\bigr\rfloor+2$ in most cases, and occasionally one less,
namely, when $3^{e-1}<(f+1)/2\le 3^{e-1}+e/4$ for some integer $e$.
However, other choices of $e$ produce weaker but interesting results.

In particular, to confirm the conjectured congruence~\eqref{eq:conj}
we may set $e=\lfloor(f+1)/2\rfloor$ in Theorem~\ref{thm:main}.
This is allowed for
$f=1$ or $f>2$, but for $f=2$ as well if we take Remark~\ref{rem:p-adic} into account.
Then $f-2e+2\ge 1$,
and congruence~\eqref{eq:conj} follows
because $\beta\equiv 1\pmod{3}$
and
$\binom{2q}{q}\equiv 2\pmod{p}$ if $q$ is a power of a prime $p$.

By applying Theorem~\ref{thm:main} with $e=\lfloor(f-3)/2\rfloor$
and using a congruence of Jacobsthal
we prove the following refinement of congruence~\eqref{eq:conj}.
\begin{cor}\label{cor:a=2}
If $q>3$ is a power of $3$ we have
\[
\frac{1}{q^2}\sum_{k=0}^{q-1}\binom{2k}{k}
\equiv
217\pmod{3^5}.
\]
\end{cor}

Because $217=1+8\cdot 3^3$, Corollary~\ref{cor:a=2} implies that congruence~\eqref{eq:conj} actually holds modulo $3^3$.
Note that the congruence of Corollary~\ref{cor:a=2} becomes an equality when $q=9$, while congruence~\eqref{eq:conj} becomes an equality when $q=3$.
We summarize further refinements in the following statement, which is obtained by applying Theorem~\ref{thm:main} with $e=\lfloor(f+1)/2\rfloor-a$.

\begin{cor}\label{cor:more}
If $a$ is a positive integer then
\[
\frac{1}{q^2}\sum_{k=0}^{q-1}\binom{2k}{k}
=-\beta\gamma
\pmod{3^{2a+1}}
\]
for $f\ge 2a+2\log(4a)/\log 3+2$,
where $\gamma$ is the $3$-adic integer
$\displaystyle
\gamma=\lim_{f\to\infty}\binom{2\cdot 3^f}{3^f}.
$
\end{cor}

Machine calculations and heuristics suggest that the congruence of Corollary~\ref{cor:more} should hold for $f\ge a$, but not for $f\ge a-1$.
(This justifies our choice of the modulus $3^{2a+1}$.)
However, our method falls short of proving this, by a finite amount of computation for each value of $a$.
A machine calculation shows that
\[
-\beta\gamma=
1+2\cdot 3^3+2\cdot 3^4+2\cdot 3^5+3^6+3^7+3^{10}+\cdots.
\]

In Section~\ref{sec:main} we present the main arguments which prove the results stated here.
To avoid distracting the reader we postpone to Section~\ref{sec:auxiliary} the statements
and proofs of some auxiliary results.
These include elementary congruences for binomial coefficients which we need in the course of the proofs.
A less elementary ingredient is Lemma~\ref{lemma:Lehmer},
whose proof depends on the use of generating functions in a $p$-adic context.

\section{Proofs of the main results}\label{sec:main}

Our proof of Theorem~\ref{thm:main} starts with some direct manipulations of the sum $\sum_{k=0}^{n-1}\binom{2k}{k}$.

\begin{proof}[Proof of Theorem~\ref{thm:main}]
Evaluating the polynomial identity
\[
\sum_{k\le m}\binom{m+r}{k}x^ky^{m-k}=
\sum_{k\le m}\binom{-r}{k}(-x)^k(x+y)^{m-k},
\]
which is~\cite[Equation~(5.19)]{GKP},
on $x=4$ and $y=-3$ and reading it from right to left we obtain
\begin{align*}
\sum_{k=0}^{n-1}\binom{2k}{k}
&=
\sum_{k=0}^{n-1}(-4)^k\binom{-1/2}{k}
\\&=
\sum_{k=0}^{n-1}4^k(-3)^{n-1-k}\binom{n-1/2}{k}
\\&=
\sum_{k=1}^{n}4^{n-k}(-3)^{k-1}\binom{n-1/2}{n-k}.
\end{align*}
The advantage of this last sum over the original sum is that, because of the power of $3$ involved,
its $3$-adic order and leading term are determined by its initial terms.
A standard identity, which can be found in~\cite[Equation~(5.36)]{GKP},
shows that the very first term of the sum, where $k=1$, equals
\[
4^{n-1}\binom{n-1/2}{n-1}
=2^{2n-1}n\binom{n-1/2}{n}
=\frac{n}{2}\binom{2n}{n}.
\]
We collect it from our sum using
\begin{align*}
\binom{n-1/2}{n-k}
&=
2^{k-1}\frac{(n-1)(n-2)\cdots(n-k+1)}{3\cdot 5\cdots(2k-1)}\binom{n-1/2}{n-1}
\\&=
2^{k-1}\frac{(n-1)!}{(n-k)!}\cdot\frac{2^k\,k!}{(2k)!}\binom{n-1/2}{n-1}
\\&=
\frac{2^{2k-1}}{k}\binom{n-1}{k-1}\binom{2k}{k}^{-1}\binom{n-1/2}{n-1},
\end{align*}
and hence obtain
\begin{equation}\label{eq:identity}
\sum_{k=0}^{n-1}\binom{2k}{k}
=
n\binom{2n}{n}\sum_{k=1}^{n}\frac{(-3)^{k-1}}{k}
\binom{2k}{k}^{-1}\binom{n-1}{k-1}.
\end{equation}

Now set $n=q$, a power of $p=3$.
Because
$(-3)^{k-1}k^{-1}\binom{2k}{k}^{-1}$
is a $3$-adic integer for all $k\ge 1$,
as follows from Lemma~\ref{lemma:p-adic}, for example,
Equation~\eqref{eq:identity} already shows that
$\sum_{k=0}^{q-1}\binom{2k}{k}$ is a multiple of $q$.
To proceed further we need to consider $3$-adic estimates of the individual terms of the sum
at the right-hand side of Equation~\ref{eq:identity}.
For that purpose we split the summation range into two segments $0<k\le 4\cdot 3^{e-1}$ and $4\cdot 3^{e-1}<k\le q$,
where $e$ is a positive integer, hence satisfying $e<f$.
The choice of the splitting point is motivated by Lemma~\ref{lemma: binom_partial}.

We start by noting that the restriction of the sum to the higher portion of the range vanishes modulo a high power of $3$.
In fact, if $e>1$ then Lemma~\ref{lemma:p-adic} implies
\begin{equation}\label{eq:higher}
\frac{3^{k-1}}{k}
\binom{2k}{k}^{-1}
\equiv
0
\pmod{3^{4\cdot 3^{e-1}-e}}
\qquad\text{for $k>4\cdot 3^{e-1}$}.
\end{equation}
The exponent of $3$ in the modulus is due the fact that
$3^e$ is the largest power of $3$ not exceeding
$2(4\cdot 3^{e-1}+1)-1$ if $e>1$.
This shows that the congruence holds for $k=4\cdot 3^{e-1}+1$.
Beyond that value of $k$ we can proceed by induction:
each time $k$ is incremented by one, the new factor $3^{k-1}$ would actually allow us to multiply the modulus by $3$,
possibly except when the new $2k-1$ equals a power of $3$.
When $e=1$, congruence~\eqref{eq:higher} only holds modulo $3^2$,
instead of $3^3$ as the formula would read, but we will see in the next paragraph that this difference is immaterial.

Now we deal with the lower part of the range.
According to Lemma~\ref{lemma: binom_partial}, we have
\[
\binom{q}{k}
\equiv (-1)^{k-1}q/k\pmod{3^{2f-2e+2}}
\]
for $0<k<4\cdot 3^{e-1}$,
and hence
\begin{equation}\label{eq:binomial}
\binom{q-1}{k-1}
=\sum_{i=0}^{k-1}(-1)^{k+i-1}\binom{q}{i}
\equiv
(-1)^{k-1}(1-qH_{k-1})
\pmod{3^{2f-2e+2}}
\end{equation}
for $0<k\le 4\cdot 3^{e-1}$.
To be able to neglect the higher portion of the summation range based on congruence~\eqref{eq:higher},
we now assume that $4\cdot 3^{e-1}-e\ge 2f-2e+2$, which amounts to the hypothesis stated in Theorem~\ref{thm:main}.
Note that when $e=1$ we actually need to assume the strict inequality $5>2f+2$, but this makes no difference
as $f$ is a positive integer.

Under this assumption, and noting that $qH_{k-1}\equiv 0\pmod{p}$ for $k\le q$,
Equations~\eqref{eq:identity},~\eqref{eq:higher} and~\eqref{eq:binomial} yield
\begin{multline*}
q^{-1}\binom{2q}{q}^{-1}
\sum_{k=0}^{q-1}\binom{2k}{k}
=
\sum_{k=1}^{q}\frac{(-3)^{k-1}}{k}
\binom{2k}{k}^{-1}\binom{q-1}{k-1}
\\\equiv
\sum_{k=1}^{m}\frac{3^{k-1}}{k}
\binom{2k}{k}^{-1}
-q\sum_{k=2}^{m}\frac{3^{k-1}}{k}
\binom{2k}{k}^{-1}H_{k-1}
\pmod{3^{2f-2e+2}},
\end{multline*}
where $m$ is any integer with $m\ge 4\cdot 3^{e-1}$.
The crucial Lemma~\ref{lemma:Lehmer} shows that the first sum in the last expression is congruent to zero,
and the desired conclusion follows.
\end{proof}

\begin{rem}\label{rem:p-adic}
The modulus in congruence~\eqref{eq:higher} is best possible when $e=2$, but can be increased to $3^{4\cdot 3^{e-1}}$ when $e>2$.
This can be proved by adapting the proof of Lemma~\ref{lemma:p-adic} to the expression being evaluated in congruence~\eqref{eq:higher},
in the specific range for $k$ under consideration.
However, it would only give a marginal improvement to Theorem~\ref{thm:main} for $q$ large
(relaxing the condition on $e$ to  $4\cdot 3^{e-1}+2e\ge 2f+2$ for $e>2$), at the expense of a tedious argument.
\end{rem}

The proof of Theorem~\ref{thm:main} already shows that the congruence in Proposition~\ref{prop:main}
holds for any $s\ge 4\cdot 3^{e-1}$.
This bound is nearly optimal if $e$ is chosen as small as possible,
but can be lowered otherwise, as we prove now.

\begin{proof}[Proof of Proposition~\ref{prop:main}]
As we did for congruence~\eqref{eq:higher} in the Proof of Theorem~\ref{thm:main},
we prove by induction on $k$ that
\begin{equation}\label{eq:sufficient}
\frac{3^{k-1}}{k}
\binom{2k}{k}^{-1}H_{k-1}
\equiv 0
\pmod{3^{f-2e+2}}
\qquad\text{for $k>s$},
\end{equation}
where $s$ satisfies an appropriate condition, to be determined in the course of the proof.

According to Lemma~\ref{lemma:p-adic}, the factor
$\frac{3^{k-1}}{k}
\binom{2k}{k}^{-1}$
is congruent to zero modulo $3^{k-1}$ divided by the largest power of $3$ not exceeding $2k-1$.
This remains true after multiplying by the other factor $H_{k-1}$ provided we further divide the modulus by the largest power of $3$ not exceeding $k-1$.
Hence the case $k=s+1$ of congruence~\eqref{eq:sufficient} holds if
$s-\lfloor\log(2s^2+s)/\log(3)\rfloor\ge f-2e+2$.
Beyond that initial value of $k$, each time $k$ is incremented by one we can actually multiply by $3$ the modulus in congruence~\eqref{eq:sufficient},
possibly except when either $k-1$ or $2k-1$ equals a power of $3$ (for the new value of $k$),
which will never occur simultaneously.

We can slightly weaken and simplify the condition found for $s$ by arguing as follows.
Taking into account the second statement in Lemma~\ref{lemma:p-adic}, the above argument shows that
if $2s+1$ is not a power of $3$ then
a sufficient condition for congruence~\eqref{eq:sufficient} to hold is
$s-\lfloor\log(2s^2)/\log 3\rfloor\ge f-2e+2$.
However, one easily sees that if $2s+1$ is a power of $3$ this condition on $s$ is equivalent to the condition found earlier.
Finally, because $2s^2$ cannot be a power of $3$ our condition is equivalent to
$s-\log(2s^2)/\log 3\ge f-2e+1$,
that is,
$3^s/(2s^2)\ge 3^{f-2e+1}$.
\end{proof}

\begin{proof}[Proofs of Corollaries~~\ref{cor:a=2} and~\ref{cor:more}]
We start with proving the more general Corollary~\ref{cor:more}.
We set $e=\lfloor(f+1)/2\rfloor-a$ in Theorem~\ref{thm:main}, so that $f-2e+2$ equals $2a+1$ when $f$ is odd, and $2a+2$ when $f$ is even.
Then Theorem~\ref{thm:main} can be applied as soon as
$4\cdot 3^{e-1}\ge 3e+4a$ when $f$ is odd, and
$4\cdot 3^{e-1}\ge 3e+4a+2$ when $f$ is even.
These are satisfied for $f$ sufficiently large,
a rough sufficient condition being
$f\ge 2a+2\log(4a)/\log 3+2$.
Now Lemma~\ref{lemma:Jacobsthal} shows that
the $3$-adic limit $\gamma$ of $\binom{2q}{q}$ exists, for $q=3^f$ tending to infinity, and also that
$\binom{2q}{q}\equiv\gamma\pmod{3^{2a+1}}$
if $f$ satisfies the above condition, because then $3f+2\ge 2a+1$.
This completes the proof of Corollary~\ref{cor:more}.
For even $f$, in the stated range, we have actually proved the congruence modulo $3^{2a+2}$.

To prove Corollary~\ref{cor:a=2} we set $a=2$.
Because $-\beta\gamma\equiv 217\pmod{3^5}$, the desired conclusion follows from Corollary~\ref{cor:more} for $f$ sufficiently large,
namely, for $f\ge 9$ if we use the more precise sufficient conditions stated above.
For the remaining values $1<f<9$ the congruence can be and has been checked by computer.
This is a matter of a few seconds for $f\le 7$, and a few minutes for $f=8$.
However, the necessary calculations can also be speeded up by an appropriate use of Equations~\eqref{eq:identity} and~\eqref{eq:higher}.
\end{proof}

\section{Some auxiliary results}\label{sec:auxiliary}

\begin{lemma}\label{lemma:p-adic}
Let $p$ be a prime and $k$ a positive integer.
Then the highest power of $p$ which divides
$k\binom{2k}{k}$
does not exceed
$2k-1$ if $p>2$, and $k+1$ if $p=2$.
Equality in these bounds is attained exactly when $2k-1$ is a power of $p>2$, or $k+1$ is a power of $p=2$.
\end{lemma}

\begin{proof}
According to a well-known theorem of Kummer, the $p$-adic order $\ord_p\bigl(\binom{n}{k}\bigr)$
of a binomial coefficient $\binom{n}{k}$, with $0\le k\le n$, equals the number
of carries which occur in adding up $k$ and $n-k$ written in base $p$.
The smallest $k$ which requires $a$ carries to add to itself in base $p$ is $k=(p^a-1)/2+1$ if $p$ is odd,
and $k=2^a-1$ if $p=2$.
For $p>2$ it follows that the largest power of $p$ dividing
$\binom{2k}{k}$ does not exceed $2k-1$.
The same holds for $k\binom{2k}{k}$ if $k$ is prime to $p$.
The general case follows from this by noting that if $k$ is a multiple of $p$ then
$k\binom{2k}{k}$ has the same $p$-adic order as
$p\cdot(k/p)\binom{2k/p}{k/p}$, again according to Kummer's theorem.
Similar arguments apply for $p=2$.
The assertion about equality in the bounds follows as well.
\end{proof}

If $q$ is a power of a prime $p$, it is well known and easy to see that
$
\binom{q}{k}
\equiv (-1)^{k-1}q/k\pmod{p^2}
$
for $0<k<q$.
This actually holds modulo $8$ when $p=2$, but for odd $p$ one can prove that the modulus is best possible.
However, the modulus can be increased if one only requires the congruence to hold on a shorter range, as in the following result.

\begin{lemma}\label{lemma: binom_partial}
Let $p^e\le p^f=q$ be powers of a prime $p$.
Then
\[
\binom{q}{k}
\equiv (-1)^{k-1}q/k\pmod{p^{2f-2e+2}}
\]
for $0<k<p^e$.

Furthermore, if $p>2$ and $e<f$ the congruence holds
for $0<k<p^e+p^{e-1}$.
\end{lemma}

\begin{proof}
We have
\begin{equation}\label{eq:A}
\binom{q}{k}=\frac{q}{k}\binom{q-1}{k-1}
\equiv 0\pmod{p^{f-e+1}}
\end{equation}
for $0<k<p^e$.
It follows that
\begin{equation}\label{eq:B}
\binom{q-1}{k-1}=\sum_{i=0}^{k-1}(-1)^{k+i-1}\binom{q}{i}
\equiv (-1)^{k-1}\pmod{p^{f-e+1}}
\end{equation}
for $0<k\le p^e$,
and hence
\begin{equation}\label{eq:C}
\binom{q}{k}=\frac{q}{k}\binom{q-1}{k-1}
\equiv (-1)^{k-1}q/k\pmod{p^{2f-2e+2}}
\end{equation}
for $0<k<p^e$.

Now assume $p>2$ and $e<f$.
Then $\sum_{i=1}^{p-1}1/i\equiv 0\pmod{p}$,
and hence
$\sum_{i=1}^{p^e-1}q/i\equiv 0\pmod{p^{f-e+2}}$.
Using congruence~\eqref{eq:C} and this fact we can refine the case $k=p^e$ of congruence~\eqref{eq:B} to
\[
\binom{q-1}{p^e-1}=\sum_{i=0}^{p^e-1}(-1)^{i}\binom{q}{i}
\equiv
1-\sum_{i=0}^{p^e-1}\frac{q}{i}
\equiv 1\pmod{p^{f-e+2}},
\]
which in turn allows us to extend congruence~\eqref{eq:C} to $k=p^e$.
In particular, we have
$\binom{q}{p^e}\equiv q/p^e\pmod{p^{f-e+1}}$.
Because congruence~\eqref{eq:A} holds also for $p^e<k<2p^e$, we have
\begin{equation*}
\binom{q-1}{k-1}=\sum_{i=0}^{k-1}(-1)^{k+i-1}\binom{q}{i}
\equiv (-1)^{k-1}(1-q/p^e)\pmod{p^{f-e+1}}
\end{equation*}
for $p^e<k\le 2p^e$.
This allows us to extend the range of congruence~\eqref{eq:C} to
$0<k<p^e+p^{e-1}$, as desired.
\end{proof}

\begin{rem}
By continuing the above argument one sees that the congruence of Lemma~\ref{lemma: binom_partial} does not extend to $k=p^e+p^{e-1}$,
and hence the longer range given for $p>2$ is optimal.

When $p=2$ the argument shows that the congruence fails to hold beyond the shorter range.
However, over that range its modulus can be slightly increased, namely,
\begin{equation}\label{eq:D}
\binom{q}{k}
\equiv (-1)^{k-1}q/k\pmod{2^{2f-2e+3}}
\end{equation}
for $0<k<2^e$.
In fact, this follows from congruence~\eqref{eq:B} with the only exception of the case $k=2^{e-1}$.
However, in the range $0<k<2^e$ congruence~\eqref{eq:A} holds modulo $2^{f-e+2}$, and hence so does congruence~\eqref{eq:B},
which implies the truth of congruence~\eqref{eq:D} for $k=2^{e-1}$ as well.
\end{rem}

In this paper we only need the case $p=3$ of the following result, but it takes no more effort to deal with an arbitrary odd prime $p$.
In the following statement $\sqrt{-p}$ is a root of $x^2+p=0$ in an extension of $\Q_p$.

\begin{lemma}\label{lemma:Lehmer}
Let $p$ be an odd prime.
In the $p$-adic field $\Q_p(\sqrt{-p})$ we have
\[
\sum_{k=1}^{\infty}\frac{p^{k}}{k}\binom{2k}{k}^{-1}
=
\frac{-\sqrt{-p}}{\sqrt{1-p/4}}
\,\log\left(\sqrt{1-p/4}+\sqrt{-p}/2\right).
\]
In particular, in $\Q_3$ we have
\[
\sum_{k=1}^{\infty}\frac{3^{k}}{k}
\binom{2k}{k}^{-1}
=0.
\]
\end{lemma}

\begin{proof}
According to~\cite{Lehmer:central_binomial}, for real $x$ with $|x|<1$ we have
\begin{equation}\label{eq:Lehmer}
\frac{2\,x\,\arcsin x}{\sqrt{1-x^2}}
=
\sum_{k=1}^{\infty}\frac{(2x)^{2k}}{k}\binom{2k}{k}^{-1}.
\end{equation}
Hence this is also an identity of formal power series in $\C[[x]]$,
which we prefer to write in the equivalent form
\begin{equation}\label{eq:Lehmer_equivalent}
\frac{-2\,x}{\sqrt{1+x^2}}
\log\bigl(\sqrt{1+x^2}+x\bigr)
=
\sum_{k=1}^{\infty}\frac{(-4x^2)^{k}}{k}\binom{2k}{k}^{-1},
\end{equation}
obtained by replacing $x$ with $ix$ and using the fact that
$
-i\arcsin(ix)
=\arcsinh x
=\log\bigl(x+\sqrt{1+x^2}\bigr).
$
Of course, in the formal setting, the logarithmic expression at the left-hand side stands for the result $f\circ g$ of composing the power series without constant term
\[
g(x)=-1+\sqrt{1+x^2}+x=
x+\frac{x^2}{2}\sum_{k=0}^{\infty}\frac{(-x^2/4)^k}{k+1}\binom{2k}{k}
\]
with the logarithmic series $f(x)=\log(1+x)=-\sum_{k=1}^{\infty}(-x)^k/k$.

Because of our assumption that $p$ is odd, the series $g(x)$ has radius of convergence $r_g=1$.
Hence it converges for $x=\sqrt{-p}/2$, as $|\sqrt{-p}/2|_p=p^{-1/2}<1$.
Furthermore, the condition $M_{p^{-1/2}}(g)<r_f=1$
of~\cite[Theorem~6.1.5]{Robert} is satisfied:
writing $g(x)=\sum_{n>0}a_nx^n$ we have $M_{p^{-1/2}}(g)=\max_{n>0}\bigl(|a_n|_p\cdot p^{-n/2}\bigr)=p^{-1/2}<1$.
Consequently, we have
\[
(f\circ g)(\sqrt{-p}/2)=f\bigl(g(\sqrt{-p}/2)\bigr)
=\log\left(\sqrt{1-p/4}+\sqrt{-p}/2\right),
\]
and the claimed formula follows.

When $p=3$ the argument of the logarithm is $(1+\sqrt{-3})/2$, a cube root of unity.
The conclusion follows
because the $p$-adic logarithm vanishes on all $p^a$-th roots of unity.
\end{proof}

\begin{rem}
The formula given in Lemma~\ref{lemma:Lehmer} is meaningful also when $p=2$.
In fact, one can prove by a different method that the formula evaluates correctly in the dyadic field $\Q_2(\sqrt{-2})$,
and hence
$\sum_{k=1}^{\infty}2^k\,k^{-1}
\binom{2k}{k}^{-1}=0$
in $\Q_2$,
because the argument of the logarithm is $\sqrt{-1/2}+\sqrt{1/2}$, an eight root of unity.
However, our proof cannot cover that case, because
the component series $\sqrt{1+x^2}$ and $\log\bigl(\sqrt{1+x^2}+x\bigr)$ at the left-hand side of Equation~\eqref{eq:Lehmer_equivalent}
(as well as the series $\arcsin(x)$ at the left-hand side of Equation~\eqref{eq:Lehmer}) have convergence radius $1/2$,
and hence do not converge for $x=\sqrt{-2}/2$.
The series at the right-hand side of Equation~\eqref{eq:Lehmer_equivalent}
has radius of convergence $2$ in this case.

In this connection it may be worth noting that
$\sum_{k=1}^{\infty}4^k\,k^{-1}
\binom{2k}{k}^{-1}=-2$
in $\Q_2$.
This is still not accessible to the argument in the proof of Lemma~\ref{lemma:Lehmer},
but follows from the fact that the partial sums of this special series admit a closed form, due to the identity\[
\sum_{k=1}^{n}\frac{4^k}{2k}
\binom{2k}{k}^{-1}
=
4^n\binom{2n}{n}^{-1}-1.
\]
This identity is~\cite[Equation~(2.9)]{Gou},
and can be easily proved by induction.
\end{rem}

The final auxiliary result which we have used in Section~\ref{sec:main} is essentially due to Jacobsthal, and is a generalization of Wolstenholme's congruence
$\frac{1}{2}\binom{2p}{p}=\binom{2p-1}{p-1}\equiv 1\pmod{p^3}$,
for a prime $p>3$ (but only modulo $3^2$ for $p=3$).

\begin{lemma}\label{lemma:Jacobsthal}
Let $p^e\le q=p^f$ be powers of an odd prime $p$.
If $p>3$ then
\[
\binom{2q}{q}\equiv\binom{2p^e}{p^e}
\pmod{p^{3e+3}}.
\]
If $p=3$ the congruence holds modulo $p^{3e+2}$.
\end{lemma}

\begin{proof}
The statement follows by induction on $f$, starting with $f=e$, using the following more general fact, due to Jacobsthal
(see Remark~\ref{rem:Jacobsthal}):
if $p$ is an odd prime and $0<k<n$ are integers, then
\begin{equation}\label{eq:Jacobsthal}
\binom{pn}{pk}\bigg/\binom{n}{k}\equiv 1
\pmod{r},
\end{equation}
where $r$ is the largest power of $p$ which divides $p^3nk(n-k)$ if $p>3$, and $p^2nk(n-k)$ if $p=3$.
\end{proof}

\begin{rem}\label{rem:Jacobsthal}
As mentioned in~\cite{Gra:organic}, which contains a proof, congruence~\eqref{eq:Jacobsthal} goes back to Jacobsthal in the early 1950's.
Because the original source is not easily accessible, the result was rediscovered by various authors, including Kazandzidis in the late 1960's,
to whom some later authors gave credit.
Among the latter authors are Robert and Zuber~\cite{RobertZuber} (see also~\cite[Chapter~7, Section~1.6]{Robert}),
who gave a proof based on properties of the Morita $p$-adic gamma function.
\end{rem}

\bibliography{References}

\def\cprime{$'$}
\providecommand{\bysame}{\leavevmode\hbox to3em{\hrulefill}\thinspace}
\providecommand{\MR}{\relax\ifhmode\unskip\space\fi MR }
\providecommand{\MRhref}[2]{%
  \href{http://www.ams.org/mathscinet-getitem?mr=#1}{#2}
}
\providecommand{\href}[2]{#2}
\begin{thebibliography}{GKP94}

\bibitem[GKP94]{GKP}
Ronald~E. Graham, Donald~E. Knuth, and Oren Patashnik, \emph{Concrete
  mathematics}, second ed., Addison-Wesley, New York, 1994.

\bibitem[Gou72]{Gou}
Henry~W. Gould, \emph{Combinatorial identities}, Henry W. Gould, Morgantown,
  W.Va., 1972, A standardized set of tables listing 500 binomial coefficient
  summations. \MR{MR0354401 (50 \#6879)}

\bibitem[Gra97]{Gra:organic}
Andrew Granville, \emph{Arithmetic properties of binomial coefficients. {I}.
  {B}inomial coefficients modulo prime powers}, Organic mathematics (Burnaby,
  BC, 1995), CMS Conf. Proc., vol.~20, Amer. Math. Soc., Providence, RI, 1997,
  pp.~253--276. \MR{MR1483922 (99h:11016)}

\bibitem[Leh85]{Lehmer:central_binomial}
D.~H. Lehmer, \emph{Interesting series involving the central binomial
  coefficient}, Amer. Math. Monthly \textbf{92} (1985), no.~7, 449--457.
  \MR{MR801217 (87c:40002)}

\bibitem[PS06]{PanSun}
Hao Pan and Zhi-Wei Sun, \emph{A combinatorial identity with application to
  {C}atalan numbers}, Discrete Math. \textbf{306} (2006), no.~16, 1921--1940.
  \MR{MR2251572 (2007d:05018)}

\bibitem[Rob00]{Robert}
Alain~M. Robert, \emph{A course in {$p$}-adic analysis}, Graduate Texts in
  Mathematics, vol. 198, Springer-Verlag, New York, 2000. \MR{MR1760253
  (2001g:11182)}

\bibitem[RZ95]{RobertZuber}
Alain Robert and Maxime Zuber, \emph{The {K}azandzidis supercongruences. {A}
  simple proof and an application}, Rend. Sem. Mat. Univ. Padova \textbf{94}
  (1995), 235--243. \MR{MR1370914 (96m:11014)}

\bibitem[ST]{SunTau:Catalan}
Zhi-Wei Sun and Roberto Tauraso, \emph{On some new congruences for binomial
  coefficients}, {\sf arXiv:math.NT/0709.1665}, to appear in Acta Arith.

\end{thebibliography}

\end{document}